\documentclass{elsart3-1-ppt}

\usepackage{amssymb,amsmath,latexsym}
\usepackage[english,francais]{babel}

\newtheorem{theorem}{Theorem}[section]

\newtheorem{e-proposition}[theorem]{Proposition}
\newtheorem{corollary}[theorem]{Corollary}
\newtheorem{e-definition}[theorem]{Definition\rm}

\newtheorem{theoreme}{Th\'eor\`eme}

\renewcommand{\theequation}{\arabic{equation}}
\def\og{\leavevmode\raise.3ex\hbox{$\scriptscriptstyle\langle\!\langle$~}}
\def\fg{\leavevmode\raise.3ex\hbox{~$\!\scriptscriptstyle\,\rangle\!\rangle$}}

\newcommand{\be}[1]{\begin{equation}\label{#1}}
\newcommand{\ee}{\end{equation}}

\newcommand{\R}{{\mathbb R}}
\renewcommand{\(}{\left(}
\renewcommand{\)}{\right)}
\newcommand{\nrm}[2]{\|#1\|_{L^{#2}(\R^2)}}

\providecommand{\abs}[1]{\left\vert #1 \right\vert}

\providecommand{\parcuad}[1]{\left[ #1 \right]}

\providecommand{\dr}[2]{\frac{\partial #1}{\partial #2}}

\newcommand{\beq}{\begin{eqnarray}}
\newcommand{\eeq}{\end{eqnarray}}
\newcommand{\ben}{\begin{eqnarray*}}
\newcommand{\een}{\end{eqnarray*}}
\newcommand{\ix}[1]{\int_{\R^2}{#1}\;dx}

\newcommand{\irdmu}[1]{\int_{\R^2}{#1}\;d\mu}
\newcommand{\irdmuM}[1]{\int_{\R^2}{#1}\;d\mu_M}
\newcommand{\ninf}{n_M}
\newcommand{\cinf}{c_M}
\def\L{{\mathcal L}}

\journal{the Acad\'emie des sciences}
\begin{document}

\centerline{}
\begin{frontmatter}

\selectlanguage{english}


\title{A functional framework for the Keller-Segel system: logarithmic Hardy-Littlewood-Sobolev and related spectral gap inequalities}

\selectlanguage{english}
\author[Ceremade]{Jean Dolbeault}
\ead{dolbeaul@ceremade.dauphine.fr}
\and
\author[Ceremade,DIM]{Juan Campos}
\ead{campos@ceremade.dauphine.fr, juanfcampos@gmail.com}

\address[Ceremade]{Ceremade (UMR CNRS no. 7534), Universit\'e Paris-Dauphine, Place de Lattre de Tassigny, 75775 Paris 16, France}
\address[DIM]{Departamento de Ingenier\'{\i}a Matem\'atica and CMM, Universidad de Chile, Casilla 170 Correo 3, Santiago, Chile}

\begin{abstract}
 \selectlanguage{english} This note is devoted to several inequalities deduced from a special form of the logarithmic Hardy-Littlewood-Sobolev, which is well adapted to the characterization of stationary solutions of a Keller-Segel system written in self-similar variables, in case of a subcritical mass. For the corresponding evolution problem, such functional inequalities play an important role for identifying the rate of convergence of the solutions towards the stationary solution with same mass.
\vskip 0.5\baselineskip\noindent
{\bf Keywords.} Keller-Segel model; chemotaxis; large time asymptotics; subcritical mass; self-similar solutions; relative entropy; free energy; Lyapunov functional; spectral gap; logarithmic Hardy-Littlewood-Sobolev inequality; Onofri's inequality; Legendre duality; best constants ---
MSC (2010): Primary: 26D10; 92C17. Secondary: 35B40

\selectlanguage{francais}
\vskip 0.5\baselineskip \noindent {\bf Un cadre fonctionnel pour le syst\`eme de Keller-Segel: in\'egalit\'e logarithmique de Hardy-Little\-wood-Sobolev et in\'egalit\'es de trou spectral reli\'ees}

\vskip 0.5\baselineskip\noindent{\bf R\'esum\'e} Cette note est consacr\'ee \`a plusieurs in\'egalit\'es fonctionnelles d\'eduites d'une forme particuli\`ere de l'in\'egalit\'e logarithmique de Hardy-Littlewood-Sobolev, qui est bien adapt\'ee \`a la caract\'erisation des solutions stationnaires d'un syst\`eme de Keller-Segel \'ecrit en variables auto-similaires, dans le cas d'une masse sous-critique. Pour le probl\`eme d'\'evolution correspondant, ces in\'egalit\'es fonctionnelles jouent un r\^ole important dans l'identification des taux de convergence des solutions vers la solution stationnaire de m\^eme masse.
\end{abstract}



\end{frontmatter}\vspace*{-1.5cm}

\selectlanguage{francais}
\section*{Version fran\c{c}aise abr\'eg\'ee}

Dans $\R^2$, l'in\'egalit\'e logarithmique de Hardy-Littlewood-Sobolev a \'et\'e \'etablie avec des constantes optimales dans \cite{MR1124215,MR1230930}. On peut l'\'ecrire sous la forme
\[
\ix{n\,\log\(\frac n{M\,\mu}\)}+\frac 2M\iint_{\R^2\times\R^2}(n(x)-M\,\mu(x))\,\log|x-y|\,(n(y)-M\,\mu(y))\;dx\,dy\ge0
\]
o\`u $M=\ix n$ et $1/\mu(x)=\pi\,(1+|x|^2)^2$ pour tout $x\in\R^2$. De plus, par dualit\'e de Legendre, elle est \'equivalente \`a l'in\'egalit\'e d'Onofri euclidienne (voir \cite{MR2433703,1101} et \cite{MR677001} pour une forme \'equivalente sur la sph\`ere).

Pour \'etudier le syst\`eme parabolique-elliptique de Keller-Segel \'ecrit en variables auto-similaires
\be{RKS-F}
\dr{n}{t} = \Delta n+\nabla\cdot(n\,x)-\nabla\cdot(n\,\nabla c)\;,\quad c=(-\Delta)^{-1}n\;,\quad x\in\R^2\,,\quad t>0\;,
\ee
on est amen\'e \`a consid\'erer une forme de l'in\'egalit\'e logarithmique de Hardy-Littlewood-Sobolev qui s'\'ecrit, sous r\'eserve que $M<8\,\pi$, sous la forme
\be{Ineq:logHLS-F}
\ix{n\,\log\(\frac n{n_M}\)}+\frac 1{4\,\pi}\iint_{\R^2\times\R^2}(n(x)-n_M(x))\,\log|x-y|\,(n(y)-n_M(y))\;dx\,dy\ge0
\ee
et o\`u $(\ninf,\cinf)$ est l'unique solution stationnaire, r\'eguli\`ere, \`a sym\'etrie radiale, de \eqref{RKS-F}, donn\'ee par
\[\label{StatKS-F}
-\Delta c_M=M\,\frac{e^{-\frac 12\,|x|^2+c_M}}{\ix{e^{-\frac 12\,|x|^2+c}}}=:n_M\;,\quad x\in\R^2\;.
\]
Exactement comme dans \cite{MR1124215,MR1230930,MR2433703,1101}, on montre par dualit\'e de Legendre qu'\`a \eqref{Ineq:logHLS-F} correspond une nouvelle in\'egalit\'e de type Onofri.
\smallskip\begin{theoreme}\label{Thm:NewOnofri-F}
Pour tout $M\in(0,8\,\pi)$, pour toute fonction $\phi$ r\'eguli\`ere \`a support compact, on a 
\[\label{Ineq:Onofri-F}
\log\(\int_{\R^2}e^\phi\,d\mu_M\)-\int_{\R^2}\phi\;d\mu_M\le\frac 1{2\,M}\ix{|\nabla\phi|^2}\;.
\]
\end{theoreme}\smallskip
Ici, $d\mu_M:=\frac 1M\,n_M\,dx$ est une mesure de probabilit\'e et comme dans \cite{DD2011}, on montre une in\'egalit\'e de trou spectral en effectuant un d\'eveloppement autour de $\phi\equiv 1$. Par densit\'e, il est par ailleurs possible d'\'etendre l'in\'egalit\'e \`a l'espace fonctionnel obtenu par compl\'etion, pour la norme $\|\phi\|^2=\ix{|\nabla\phi|^2}+(\int_{\R^2}\phi\;d\mu_M)^2$, de l'ensemble des fonctions r\'eguli\`eres \`a support compact.

Dans sa forme lin\'earis\'ee, le syst\`eme de Keller-Segel s'\'ecrit
\be{Eq:LinKS-F}
\dr{f}{t}=\frac1{\ninf}\nabla\cdot\big[\ninf\nabla(f-g\,\cinf)\big]=:\L\,f\quad\mbox{o\`u}\quad g\,\cinf=(-\Delta)^{-1}(f\,\ninf)\;.
\ee
On montre que le noyau de $\L$ est engendr\'e par une fonction $f_{0,0}$ d\'etermin\'ee par $-\Delta f_{0,0}=f_{0,0}\,\ninf$. En effectuant un d\'eveloppement limit\'e \`a l'ordre deux autour de $\ninf$, il est ais\'e de voir que
\[
\mathsf Q_1[f]:=\irdmuM{|f|^2}+\frac 1{2\,\pi}\iint_{\R^2\times\R^2}f(x)\,\log|x-y|\,f(y)\;d\mu_M(x)\,d\mu_M(y)\ge 0\;.
\]
De plus $\mathsf Q_1[f]=0$ si et seulement si $f$ est proportionnelle \`a $f_{0,0}$. On montre alors le r\'esultat suivant.
\smallskip\begin{theoreme}\label{Teo:LimHLS-F} Il existe $\kappa>1$ tel que, pour tout $f\in L^2(\R^2,d\mu_M)$, si $\irdmuM{f\,f_{0,0}}=0$, alors on a
\[
\irdmuM{f^2}\le\kappa\,\mathsf Q_1[f]\;.
\]
\end{theoreme}\smallskip
Si l'on d\'efinit maintenant $\mathsf Q_2[f]:=\langle f,\,\L\,f\rangle$, on montre une derni\`ere in\'egalit\'e de trou spectral.
\smallskip\begin{theoreme}\label{Teo:Gap-F} Pour toute fonction $f\in L^2(\R^2,f\mu_M)$ v\'erifiant $\irdmuM{f\,f_{0,0}}=0$, on a
\[
\mathsf Q_1[f]\le\mathsf Q_2[f]\;.
\]
\end{theoreme}\smallskip
Il est alors facile d'en d\'eduire que si $f$ est une solution de \eqref{Eq:LinKS-F}, alors $\mathsf Q_1[f(t,\cdot]\le\mathsf Q_1[f(0,\cdot]\,e^{-2t}$ pour tout $t\ge0$. Pour une preuve d\'etaill\'ee des Th\'eor\`emes~\ref{Teo:LimHLS-F} et~\ref{Teo:Gap-F}, on renverra \`a \cite{CD2012}. Au prix d'une estimation un peu plus compliqu\'ee bas\'ee sur la formule de Duhamel, on montre que cette estimation en temps grand s'applique aussi \`a $f:=(n-\ninf)/\ninf$, o\`u $n$ est la solution de \eqref{RKS-F}.

\selectlanguage{english}
\setcounter{equation}{0}
\par\medskip\centerline{\rule{2cm}{0.2mm}}\medskip

\section{Introduction}\label{Sec:Intro}

In $\R^2$, the logarithmic Hardy-Littlewood-Sobolev has been established with optimal constants in \cite{MR1124215} (also see~\cite{MR1230930}) and can be written as
\be{Ineq:logHLSCarlenLoss}
\ix{n\,\log\(\frac nM\)}+\frac 2M\int_{\R^2\times\R^2}n(x)\,n(y)\,\log|x-y|\;dx\,dy+M\,\(1+\log\pi\)\ge 0
\ee
for any function $n\in L^1_+(\R^2)$ with $M=\ix n$. As a consequence (see \cite{MR2103197}), the \emph{free energy} functional
\[
F[n]:=\ix{n\,\log n}+\frac 12\ix{|x|^2\,n}-\frac 12\ix{n\,c}+K\quad\mbox{with}\quad c=(-\Delta)^{-1}n:=-\frac 1{2\,\pi}\,\log|\cdot|*n
\]
is bounded from below if $M\in(0,8\,\pi]$. Here $K=K(M)$ is a constant to be fixed later. We may observe that $F$ is not bounded from below if $M>8\,\pi$, for instance by considering $\lambda\mapsto F[n_\lambda]$ where $n_\lambda(x)=\lambda^2\,n(\lambda\,x)$ for some given function $n$, and by taking the limit $\lambda\to\infty$. See \cite{0728} for more details. Equality in \eqref{Ineq:logHLSCarlenLoss} is achieved by
\[
\mu(x):=\frac 1{\pi\,(1+|x|^2)^2}\quad\forall\;x\in\R^2\,,
\]
which solves $-\Delta\log \mu=8\,\pi\,\mu$ and can be inverted as $(-\Delta)^{-1}\mu=\frac 1{8\,\pi}\,\log \mu+\frac 1{8\,\pi}\,\log\pi$. 

Consider the probability measure $d\mu:=\mu\,dx$. Written in Euclidean form, Onofri's inequality (see~\cite{MR677001} for the equivalent version on the sphere)
\be{Ineq:Onofri8pi}
\log\(\irdmu{e^{\,\phi}}\)-\irdmu \phi\le \frac 1{16\,\pi}\,\ix{|\nabla \phi|^2}
\ee
plays in dimension $d=2$ the role of Sobolev's inequality in higher dimensions. The inequality holds for any smooth function with compact support and, by density, for any function $\phi$ in the space obtained by completion with respect to the norm given by: $\|\phi\|^2=\ix{|\nabla\phi|^2}+(\int_{\R^2}\phi\;d\mu)^2$. Onofri's inequality can be seen as the \emph{dual} inequality of the logarithmic Hardy-Littlewood-Sobolev, \emph{cf} \cite{MR1124215,MR1230930,MR2433703,1101}.

\medskip The rescaled parabolic-elliptic Keller-Segel system reads
\be{RKS}
\dr{n}{t} = \Delta n+\nabla\cdot(n\,x)-\nabla\cdot(n\,\nabla c)\;,\quad c=(-\Delta)^{-1}n\;,\quad x\in\R^2\,,\quad t>0
\ee
Assume that the initial datum is $n(0,\cdot) = n_0$. If $M=\ix{n_0}>8\,\pi$, solutions blow-up in finite time. If $n_0\in L^1_+\big(\R^2\,,(1+\abs{x}^2)\,dx\big)$, $n_0\abs{\log n_0}\in L^1(\R^2)$ and $M< 8\,\pi$, solutions globally exists and it has been shown in \cite[Theorem 1.2]{MR2226917} that
\[
\lim_{t\to\infty}\nrm{n(t,\cdot)-\ninf}1=0\quad\textrm{and}\quad\lim_{t\to\infty}\nrm{\nabla c(t,\cdot)-\nabla \cinf}2=0\;,
\]
where $(\ninf,\cinf)$ is the unique, smooth and radially symmetric solution of
\be{StatKS}
-\Delta c_M=M\,\frac{e^{-\frac 12\,|x|^2+c_M}}{\ix{e^{-\frac 12\,|x|^2+c}}}=:n_M\;,\quad x\in\R^2\;.
\ee
Notice that $\ninf=M\,{e^{\cinf-\abs{x}^2/2}}/{\ix{e^{\cinf-\abs{x}^2/2}}}$ with $\cinf=(-\Delta)^{-1}\ninf$. The case $M=8\,\pi$ has also been extensively studied, but is out of the scope of this note.

Ineq.~\eqref{Ineq:Onofri8pi} and the Moser-Trudinger inequality have been repeatedly used to study the Keller-Segel system in bounded domains. In the whole space case, Ineq.~\eqref{Ineq:logHLSCarlenLoss} turns out to be very convenient, at least for existence issues. Ineq.~\eqref{Ineq:Onofri8pi} and Ineq.~\eqref{Ineq:logHLSCarlenLoss} correspond to the $M=8\,\pi$ case. For $M<8\,\pi$, we will establish a new inequality of Onofri type, which is our first main result: see Theorem~\ref{Thm:NewOnofri}.

\medskip An important issue in the study of \eqref{RKS} is to characterize the rate of convergence of $n$ towards $\ninf$. See~\cite{Blanchet2010533,Calvez:2010fk}. For this purpose, it is convenient to linearize the Keller-Segel system~\eqref{RKS} by considering
\[
n(t,x)=\ninf(x)\,(1+\varepsilon\,f(t,x))\quad\textrm{and}\quad c(t,x)=\cinf(x)\,(1+\varepsilon\,g(t,x))
\]
and formally take the limit as $\varepsilon\to0$. At order $O(\varepsilon)$, $(f,g)$ solves
\be{Eq:LinKS}
\dr{f}{t}=\frac1{\ninf}\nabla\cdot\big[\ninf\nabla(f-g\,\cinf)\big]=:\L\,f\quad\mbox{and}\quad g\,\cinf=(-\Delta)^{-1}(f\,\ninf)\;.
\ee
As we shall see in Section~\ref{Sec:Gaps}, several spectral gap inequalities (see Theorems~\ref{Teo:LimHLS} and~\ref{Teo:Gap}) are related with~\eqref{Ineq:logHLSCarlenLoss} and involve the linear operator $\L$. Detailed proofs and applications to the full Keller-Segel system~\eqref{RKS} will be given in a forthcoming paper, \cite{CD2012}, whose main result is that $\nrm{n(t,\cdot)-\ninf}1=O(e^{-t})$ as $t\to\infty$.

\section{Duality and stationary solutions of the Keller-Segel model in self-similar variables}\label{Sec:Duality}

For any $M\in(0,8\,\pi)$, the function $\cinf$ given by \eqref{StatKS} can be characterized either as a minimizer of
\[
G[c]:=\frac 12\ix{n\,c}-M\,\log\(\ix{e^{-\frac 12\,|x|^2+c}}\)
\]
where $n$ and $c$ are related through the Poisson equation, $-\Delta c=n$, or in terms of $n$, seen as a minimizer of the functional $n\mapsto F[n]$. Inspired by \cite{MR1230930,MR2433703,MR1124215,1101}, we can characterized the corresponding functional inequalities and observe that they are \emph{dual} of each other. Let us give some details.

\medskip Consider the \emph{free energy} functional $n\mapsto F[n]=F_1[n]-F_2[n]$ (for an appropriate choice of the constant $K$) on the set $\mathcal X_M$ of all nonnegative integrable functions with mass $M>0$, where
\[
F_1[n]=\ix{n\,\log\(\frac n{n_M}\)}\quad\mbox{and}\quad F_2[n]=\frac 12\ix{(n-n_M)\,(-\Delta)^{-1}(n-n_M)}\;.
\]
The free energy $F$ is bounded from below by~\eqref{Ineq:logHLSCarlenLoss}. Since $n_M$ is a minimizer for $F$ and $F[n_M]=0$, we actually have the functional inequality $F_1[n]\ge F_2[n]$ for any $n\in\mathcal X_M$. This inequality can be rewritten as
\[\label{Ineq:logHLS}
\ix{n\,\log\(\frac n{n_M}\)}+\frac 1{4\,\pi}\iint_{\R^2\times\R^2}(n(x)-n_M(x))\,\log|x-y|\,(n(y)-n_M(y))\;dx\,dy\ge0
\]
for any $n\in\mathcal X_M$ with $M<8\,\pi$.

By Legendre's duality, we have: $F_1^*[\phi]\le F_2^*[\phi]$ where $F_i^*[\phi]:=\sup_{n\in\mathcal X_M}\(\ix{\phi\,n}-F_i[n]\)$, $i=1$, $2$, is defined on $L^\infty(\R^2)$. A straightforward computation shows that $F_1^*[\phi]=\ix{\phi\,n}-F_1[n]$ if and only if $\log(\frac n{n_M})=\phi-\log\(\irdmuM{e^\phi}\)+\log M$, so that
\[
F_1^*[\phi]=M\,\log\(\irdmuM{e^\phi}\)-M\,\log M\;.
\]
Here $d\mu_M$ is the probability measure
\[
d\mu_M:=\mu_M\,dx\;,\quad\mbox{with}\quad\mu_M:=\frac 1M\,n_M\;.
\]

It is clear that we can impose at no cost that $\irdmuM{\phi}=0$. It is also standard to observe that $F_2^*[\phi]=\ix{\phi\,n}-F_2[n]$ if and only if $\phi=(-\Delta)^{-1}(n-n_M)$, so that
\[
F_2^*[\phi]=\frac 12\ix{|\nabla\phi|^2}\;.
\]
Notice that $\ix{|\nabla\phi|^2}$ is well defined as $-\Delta\phi=n-n_M$ is integrable and such that $\ix{(n-n_M)}=0$. With $c_M=(-\Delta)^{-1}n_M$ and $\phi=c-c_M$, we recover that $G[\phi+c_M]$ is equal to $F_2^*[\phi]-F_1^*[\phi]$ up to a constant. Replacing $\phi$ by $\phi-\int_{\R^2}\phi\;d\mu_M$, we arrive at the following result in the space $\mathcal H_M$ obtained by completion with respect to the norm given by: $\|\phi\|^2=\ix{|\nabla\phi|^2}+(\int_{\R^2}\phi\;d\mu_M)^2$.
\smallskip\begin{theorem}\label{Thm:NewOnofri}
For any $M\in(0,8\,\pi)$, with $n_M$ defined as the unique minimizer of $F$, {\rm i.e.} the unique solution $n_M$ given by \eqref{StatKS},  and $n_M\,dx=d\mu_M$, with $c_M=(-\Delta)^{-1}\,n_M$, we have the following inequality:
\be{Ineq:Onofri}
\log\(\int_{\R^2}e^\phi\,d\mu_M\)-\int_{\R^2}\phi\;d\mu_M\le\frac 1{2\,M}\ix{|\nabla\phi|^2}\quad\forall\;\phi\in\mathcal H_M\;.
\ee
\end{theorem}\smallskip
As a consequence, if we consider the special case $\phi=1+\varepsilon\,\psi$ and consider the limit $\varepsilon\to 0$ in \eqref{Ineq:Onofri}, as in~\cite{DD2011}, we get an interesting spectral gap inequality.
\smallskip\begin{corollary}\label{Cor:NewOnofri} With the above notations, for any $\psi\in \mathcal H_M$, the following inequality holds
\[
\int_{\R^2}\left|\psi-\overline\psi\right|^2\,n_M\,dx\le\ix{|\nabla\psi|^2}\quad\mbox{where}\quad \overline\psi=\int_{\R^2}\psi\;d\mu_M\;.
\]
\end{corollary}

\section{Linearized Keller-Segel model, spectral gap inequalities and consequences}\label{Sec:Gaps}

Exactly as for Ineq.~\eqref{Ineq:Onofri}, we observe that
\[
\mathsf Q_1[f]:=\irdmuM{|f|^2}+\frac 1{2\,\pi}\iint_{\R^2\times\R^2}f(x)\,\log|x-y|\,f(y)\;d\mu_M(x)\,d\mu_M(y)=\lim_{\varepsilon\to0}\frac 1{\varepsilon ^2}\,F[\ninf(1+\varepsilon\,f)]\ge 0
\]
Notice that $\mathsf Q_1[f]=\ix{|f|^2\,\ninf}-\ix{|\nabla(g\,\cinf)|^2}$ if $\irdmuM f=0$. We also notice that $f_{0,0}:=\partial_M\log\ninf$ generates the kernel $\mathop{Ker}(\L)$ considered as an operator on $L^2(\R^2,d\mu_M)$ and the functions $f_{1,i}:=\partial_{x_i}\log\ninf$ with $i=1$, $2$ and $f_{0,1}:=x\cdot\nabla\log\ninf$ are eigenfunctions of $\L$ with eigenvalues $1$ and $2$ respectively; moreover they generate the corresponding eigenspaces (see \cite{CD2012} for details). It is remarkable that $\mathsf Q_1[f]=0$ if and only if $f\in\mathop{Ker}(\L)$ and this allows to establish a first spectral gap inequality.
\smallskip\begin{theorem}\label{Teo:LimHLS} There exists $\kappa>1$ such that
\[
\irdmuM{f^2}\le\kappa\,\mathsf Q_1[f]\quad\forall\;f\in L^2(\R^2,f\mu_M)\quad\mbox{such that}\irdmuM{f\,f_{0,0}}=0\;.
\]
\end{theorem}\smallskip

\medskip The proof of Theorem~\ref{Teo:LimHLS} relies on spectral properties of Schr\"odinger operators. See \cite{CD2012} for details. Since $\mathsf Q_1[f]=0$ if and only if $f\in\mathop{Ker}(\L)$, that is if $f$ is proportional to $f_{0,0}$, we can define the scalar product $\langle\cdot,\cdot\rangle$ induced by the quadratic form $\mathsf Q_1$ on the space $\mathcal D_M$ orthogonal of $f_{0,0}$ in $L^2(\R^2,d\mu_M)$. With this definition, we have $\mathsf Q_1[f]=\langle f,f\rangle$. On the space $\mathcal D_M$ with scalar product $\langle\cdot,\cdot\rangle$, the operator $\mathcal L$ is self-adjoint. Let
\[
\mathsf Q_2[f]:=\langle f,\,\L\,f\rangle\;.
\]
Then we have a second spectral gap inequality.
\smallskip\begin{theorem}\label{Teo:Gap} For any function $f\in\mathcal D_M$, we have
\[
\mathsf Q_1[f]\le\mathsf Q_2[f]\;.
\]
\end{theorem}\smallskip
Moreover, if $f$ is a radial function, then we have $2\,\mathsf Q_1[f]\le\mathsf Q_2[f]$. The operator $\L$ has only discrete spectrum as a consequence of Persson's lemma, or as can be shown by direct investigation using the tools of the concentration-compactness method and the Sturm-Liouville theory. By rewriting the spectral problem for $\L$ in terms of \emph{cumulated densities,} it is possible to prove that the eigenspace corresponding to the lowest non-zero eigenvalue is generated by $f_{1,i}$ with $i=1$, $2$, which completes the proof. See \cite{CD2012} for details.

\medskip As a simple consequence, if $f$ is a solution to \eqref{Eq:LinKS}, then
\[
\frac d{dt}\langle f,f\rangle=-\langle f,\,\L\,f\rangle\le-\,2\,\langle f,f\rangle\;,
\]
which shows the exponential convergence of $f$ towards $0$. The nonlinear Keller-Segel model \eqref{RKS} can be rewritten in terms of $f:=(n-\ninf)/\ninf$ and $g:=(c-\cinf)/\cinf$ as 
\[
\dr{f}{t} -\L\,f = -\frac1{\ninf}\,\nabla\cdot\parcuad{f\,\ninf(\nabla(g\,\cinf))}\;.
\]
Estimates based on Duhamel's formula allow to prove that $t\mapsto\mathsf Q_1[f(t,\cdot)]$ is bounded uniformly with respect to $t>0$ and
\[
\frac d{dt}\,\mathsf Q_1[f(t,\cdot)]\le-\,\mathsf Q_1[f(t,\cdot)]\left[2-\delta(t,\varepsilon)\(\mathsf Q_1[f(t,\cdot)])^\frac{1-\varepsilon}{2-\varepsilon}+\mathsf Q_1[f(t,\cdot)])^\frac 1{2+\varepsilon}\)\right]\,.
\]
for any $\varepsilon>0$ small enough, for some continuous $\delta$ such that $\lim_{t\to\infty}\delta(t,\varepsilon)=0$. This proves that $\lim_{t\to\infty}e^{2t}\,\mathsf Q_1[f(t,\cdot)]$ is finite. Details will be given~in~\cite{CD2012}.

\linespread{0.9}

\begin{thebibliography}{10}

\bibitem{MR1230930}
{\sc W.~Beckner}, {\em Sharp {S}obolev inequalities on the sphere and the
  {M}oser-{T}rudinger inequality}, Ann. of Math. (2), 138 (1993), pp.~213--242.

\bibitem{Blanchet2010533}
{\sc A.~Blanchet, J.~Dolbeault, M.~Escobedo, and J.~Fern{\'a}ndez}, {\em
  Asymptotic behaviour for small mass in the two-dimensional parabolic-elliptic
  {K}eller-{S}egel model}, J. Math. Analysis and Applications, 361 (2010),
  pp.~533 -- 542.

\bibitem{MR2226917}
{\sc A.~Blanchet, J.~Dolbeault, and B.~Perthame}, {\em Two-dimensional
  {K}eller-{S}egel model: optimal critical mass and qualitative properties of
  the solutions}, Electron. J. Differential Equations, 44 (2006), pp.~1--32
  (electronic).

\bibitem{Calvez:2010fk}
{\sc V.~Calvez and J.~A. Carrillo}, {\em Refined asymptotics for the
  subcritical {K}eller-{S}egel system and related functional inequalities}.
\newblock Preprint ArXiv 1007.2837, to appear in Proc. AMS.

\bibitem{MR2433703}
{\sc V.~Calvez and L.~Corrias}, {\em The parabolic-parabolic {K}eller-{S}egel
  model in {$\mathbb R\sp 2$}}, Commun. Math. Sci., 6 (2008), pp.~417--447.

\bibitem{CD2012}
{\sc J.~Campos and J.~Dolbeault}, {\em Asymptotic estimates for the
  parabolic-elliptic {K}eller-{S}egel model in the plane}.
\newblock Preprint, 2012.

\bibitem{MR1124215}
{\sc E.~A. Carlen and M.~Loss}, {\em Competing symmetries of some functionals
  arising in mathematical physics}, in Stochastic processes, physics and
  geometry ({A}scona and {L}ocarno, 1988), World Sci. Publ., Teaneck, NJ, 1990,
  pp.~277--288.

\bibitem{DD2011}
{\sc M.~{Del Pino} and J.~{Dolbeault}}, {\em {The Euclidean Onofri inequality
  in higher dimensions}}, arXiv preprint 1201.2162, to appear in Int. Math.
  Res. Notices,  (2012).

\bibitem{1101}
{\sc J.~Dolbeault}, {\em {S}obolev and {H}ardy-{L}ittlewood-{S}obolev
  inequalities: duality and fast diffusion}, Math. Research Letters,  (2011).

\bibitem{MR2103197}
{\sc J.~Dolbeault and B.~Perthame}, {\em Optimal critical mass in the
  two-dimensional {K}eller-{S}egel model in {$\mathbb R\sp 2$}}, C. R. Math.
  Acad. Sci. Paris, 339 (2004), pp.~611--616.

\bibitem{0728}
{\sc J.~Dolbeault and C.~Schmeiser}, {\em The two-dimensional {K}eller-{S}egel
  model after blow-up}, Discrete and Continuous Dynamical Systems, 25 (2009),
  pp.~109--121.

\bibitem{MR677001}
{\sc E.~Onofri}, {\em On the positivity of the effective action in a theory of
  random surfaces}, Comm. Math. Phys., 86 (1982), pp.~321--326.

\end{thebibliography}

\smallskip\noindent{\bf Acknowledgments.} The authors acknowledge support by the ANR projects \emph{CBDif-Fr} and \emph{EVOL} (JD), and by the MathAmSud project \emph{NAPDE} (JC and JD).

\noindent\copyright\ 20012 by the authors. This paper may be reproduced, in its entirety, for non-commercial purposes.

\medskip\begin{flushright}{\sl\today}\end{flushright}
\end{document}